\documentclass[a4paper,12pt]{article}
\usepackage{amssymb,hyperref}
\usepackage{color}
\newtheorem{thm}{Theorem}
\newtheorem{prop}{Proposition}
\definecolor{Red}{rgb}{1,00,0}

\begin{document}
\def\real{\mathbb{R}}
\def\nat{\mathbb{N}}
\def\fine{\hspace{\stretch{1}} $\Box$}

\begin{centering}
\LARGE{Explicit formulas using partitions of integers for numbers
defined by recursion}
 \vspace*{6mm}

 \Large{Giuseppe Fera$^\star$,\ \ Vittorino Talamini$^{\star,\ast}$}\\
 \vspace*{6mm}
 \normalsize
 {$\star$ DCFA, Sezione di Fisica e Matematica, Universit\`a di Udine, Via delle Scienze 206, 33100 Udine,
 Italy}\\
 {$\ast$ INFN, Sezione di Trieste, Via Valerio 2, 34127 Trieste, Italy}\\
 \vspace*{6mm}
 {\normalsize{\tt giuseppe.fera@uniud.it,\ \ vittorino.talamini@uniud.it}}\\
\end{centering}
\vspace*{2\baselineskip}

\begin{abstract}
In this article we obtain an explicit formula in terms of the partitions of the positive integer $n$ to express the
$n$-th term of a wide class of sequences of numbers defined by
recursion. Our proof is based only on arithmetics. 
We compare our result with similar
formulas obtained with different approaches already in the XIX
century. Examples are given for Bernoulli, Euler and Fibonacci numbers.
\end{abstract}

\section{Introduction and main result}\label{Intro}
\vskip 1cm
 Consider two sequences of numbers
$a=\{a_0,a_1,\ldots\}$ and $b=\{b_0,b_1,\ldots\}$ such that for
all integers $n\geq 0$ the following condition holds:
\begin{equation}\label{condition}
 \sum_{h=0}^{n}\,a_{n-h}\,b_h=\delta_{0,n}\,,
\end{equation}
in which $\delta_{0,n}$ is $1$ if $n=0$ and $0$ otherwise. We stress
that for $n=0$, Condition (\ref{condition}) implies $a_0 b_0=1$.
It is convenient to consider $a_0=b_0=1$. If this is not the case,
one may 
consider the sequences obtained from $a$ and $b$ by
factorizing out the factors $a_0$ and $b_0$, that is with
coefficients $a_i/a_0$ and $b_i/b_0$, $\forall i\geq 0$.
Therefore, from now on we will suppose $a_0=b_0=1$.

If, with the two given sequences, one builds the formal power
series
\begin{equation}\label{serieformali}
   a(x)=\sum_{i=0}^\infty a_i x^i,\qquad b(x)=\sum_{i=0}^\infty b_i
x^i\,,
\end{equation}
and makes the ``Cauchy product'' of the two series (with formal
power series one does not bother about convergence problems), then
Condition (\ref{condition}) is equivalent to the equation
$$a(x)\,b(x)=1\,.$$

Condition (\ref{condition}) is also equivalent to a recursive
definition of the elements 
of one of the two sequences. One has, in fact:
\begin{equation}\label{recursion}
b_n=-\sum_{h=0}^{n-1}\,a_{n-h}\,b_h\,,\qquad \forall n>0, \qquad
b_0=1\,,
\end{equation}
that allows to determine the $n$-th element of the sequence $b$ from
all the preceding ones. 
\\

Many sequences of numbers are defined through Eqs.
(\ref{condition}) or (\ref{recursion}) (some examples are given in Table \ref{tab1}),
and in those cases to calculate the $n$-th element of the sequence
one has to calculate all the elements of the sequence with index
smaller than $n$.

Formulas that allow to calculate the $n$-th element of the
sequence  only in terms of $n$, and not in terms of other elements
of the sequence, are then useful and welcome, often more for
theoretical than for practical interest. These formulas are often
called {\em closed} formulas, but we prefer to call them {\em
explicit} formulas.

In this article we present, in Theorem \ref{uno}
, an explicit formula that allows to calculate the $n$-th element of
a sequence of numbers satisfying Condition (\ref{condition}) in
terms of the partitions of the integer $n$. This is not a new
result, as we will see, but we will prove it in an elementary way
that seems absent in the literature. Later on 
we will list and comment alternative ways to get the
same results.\\

Before to state our theorem we recall some basic
definitions concerning partitions and compositions. 

Given a positive integer $n$, a {\em composition} $c$ of $n$ is an
unordered set of $l(c)$ positive
integers $n_i$ that have sum $n$:
$$c=\{n_1,n_2,\ldots,n_{l(c)}\},\qquad n=\sum_{i=1}^{l(c)} n_i\,.$$
The number $l(c)$ is called the {\em length} of the composition
$c$.

A {\em partition} $p$ is a composition in which all the $n_i$ are
ordered, usually in decreasing order. The length $l(p)$ of a
partition $p$ is the number of elements of $p$. If $p$ is a partition obtained from
the composition $c$, one obviously has $l(p)=l(c)$. If $p=\{n_1,n_2,\ldots,n_{l(p)}\}$, one then has:
$$n_i\geq n_{i+1},\qquad \forall i=1,\ldots, l(p)-1.$$
For any partition $p$ there are many compositions formed with the
same $l(p)$ elements of $p$. Let $\mu(p)$ be
the number of different compositions formed by the same numbers in
the partition $p$. It is an easy exercise to verify that the
number $\mu(p)$ is given by the formula:
\begin{equation}\label{numpart}
 \mu(p)={l(p)!}\left/{\prod_{n_i\,\in\, \cup(p)} m_p(n_i)!}\right.\,,
\end{equation}
where $m_p(n_i)$ is the {\em multiplicity} of $n_i$ in $p$, that
is the number of times the number $n_i$ appears in $p$, and the
product at the denominator has to be done for all {\em different}
numbers $n_i$ that appear in $p$, that is, for all numbers $n_i$
that appear in the union $\cup(p)$ of the elements of $p$. The
multiplicities $m_p(n_i)$ that appear at the denominator of
(\ref{numpart}) always form a composition of the length $l(p)$
that appears at the numerator. The second member of Eq.
(\ref{numpart}) can also be written using the multinomial
coefficient: $$\mu(p)={l(p) \choose {m_p(n_1),\,
m_p(n_2),\ldots}}\,,$$ where the numbers ${m_p(n_1),\,
m_p(n_2),\ldots}$ are those appearing in the denominator of the
second member of Eq. (\ref{numpart}).

Let ${\it C}(n)$  and ${\cal P}(n)$ be the set of all
compositions of $n$ and the set of all partitions of $n$,
respectively.\\

We are ready to state our main result.

\begin{thm}$\!\!\!{\bf .}$\label{uno}
Let $\{a_0,a_1,\ldots\}$ and $\{b_0,b_1,\ldots\}$ be sequences
such that $a_0=b_0=1$, and such that Condition (\ref{condition})
holds, then one may calculate the elements $b_n$, $\forall\,n>0$,
from the following explicit formulas involving the sequence
$\{a_0,a_1,\ldots\}$ and the compositions or the partitions of
$n$:
 \begin{equation}\label{eqt1}
    b_n= \sum_{c\,\in\, {\cal C}(n)} \;{\prod_{n_i\,\in\, c}
                                                \left(-{a_{n_i}}\right)}
 \end{equation}
 \begin{equation}\label{eqt2}
    b_n= \sum_{p\,\in\, {\cal P}(n)} \,\mu(p)\;{\prod_{  n_i\,\in\, p } \left(-{a_{n_i}}\right)}
 \end{equation}
\end{thm}

Formulas (\ref{eqt1}) and (\ref{eqt2}) in Theorem \ref{uno}
differ only for the use of compositions in place of partitions. In
general, the use of partitions is preferable in calculations
because it involves fewer terms. In Table \ref{tab2} we list
explicit formulas for the general terms of some sequences of
numbers in terms of partitions only, but there exist similar
formulas using compositions, too.

It must be said that the explicit formulas that use partitions of integers to calculate the general $n$-th term of a sequence of numbers are often
not convenient for practical computations, because the cardinality
of ${\cal P}(n)$ grows very rapidly with $n$.\\

\noindent {\textbf{Proof of Theorem \ref{uno} (elementary version)}}.
We want to solve the linear system expressed by Condition
(\ref{condition}) in the indeterminates $b_1,b_2,\ldots$. We
recall that we suppose $a_0=b_0=1$. Following Scherk
\cite{scherk1825}, we do this by substitution, using recursive
formula (\ref{recursion}) to solve the equation for $b_1$, then
that one for $b_2$, and so on. At each step, the equation number
$n$ contains only the indeterminates $b_1,\ldots,b_n$, so the
substitution method is easy to apply.

The first few equations are not difficult to solve:
\begin{equation}\label{sistscherk}
    \begin{array}{lcl}
     b_1&=& -a_1 \\
 b_2&=& -a_2- a_1b_1=-a_2+a_1^2\\
     b_3&=& -a_3-a_2b_1-a_1b_2=-a_3+2a_2a_1-a_1^3\\
     b_4 &=& -a_4-a_3b_1-a_2b_2-a_1b_3=-a_4+2a_3a_1+a_2^2-3a_2a_1^2+a_1^4\\
     b_5 &=& -a_5-a_4b_1-a_3b_2-a_2b_3-a_1b_4=\\
     &=&-a_5+2a_4a_1+2a_3a_2-3a_3a_1^2-3a_2^2a_1+4a_2a_1^3-a_1^5\\
      \ldots &&\\
   \end{array}
\end{equation}
At this point one notes that, at least for $h=1,2,3,4,5,\ldots$,
$b_h$ is obtained as the sum of all products like
$${a_{n_1}}\,{a_{n_2}}\,\ldots \,{a_{n_k}}\,,$$
whose indices satisfy the conditions
$$1\leq k\leq h,\qquad n_1+n_2+\ldots +n_k=h,$$
and each term is taken with the $+$ sign if $k$ is even, with the
$-$ sign if $k$ is odd. The $k$ positive numbers $n_1,\ldots,n_k$
form then a composition of length $k$ of the positive integer $h$.
We are then led to write the formula:
\begin{eqnarray}\label{enne0x}
 b_n&=&\sum_{c\,\in\, {{\cal C}(n)}}\; \prod_{n_i \,\in\, c}({-a_{n_i}}),\quad \forall n\geq 1\,,
\end{eqnarray}
where ${\cal C}(n)$ is the set of all compositions of the integer
$n$. Let's verify by induction that Formula (\ref{enne0x}) is true
for all values of $n\geq 1$. From (\ref{sistscherk}) we see that Eq. (\ref{enne0x}) is true for $n\leq 5$. Supposing (\ref{enne0x}) true up to a given $n\geq 1$, let's prove that it is true also for $n$
substituted by $n+1$. Starting from (\ref{recursion}), and using
(\ref{enne0x}), one finds:
$$b_{n+1}=-a_{n+1}-\sum_{h\,=\,1}^{n}\,a_{n+1-h}\, b_h=$$
$$=-a_{n+1}-\sum_{h\,=\,1}^{n}\,\left(\sum_{c\,\in\, {{\cal C}(h)}}\; \prod_{n_i \,\in\, c}({-a_{n_i}})\right)a_{n+1-h}.$$
For any composition $c=\{n_1,n_2,\ldots,n_{l(c)}\}\in{{\cal
C}(h)}$, consider the set obtained by appending to the $l(c)$
numbers $n_i$ in $c$ the number $n_0=(n+1)-h$, that appears in the
index of the last factor $a_{n+1-h}$ (so that it can be rewritten
in the form $a_{n_0}$). The set
$c'=\{n_0,n_1,n_2,\ldots,n_{l(c)}\}$  is now a composition of
$(n+1)$ of length $l(c')=l(c)+1$. Considering then all possible
values of $h$ resulting from the sum of the last expression, we
obtain in this way all possible compositions of $(n+1)$ of length
greater than 1. Considering the term outside the summation symbol
as corresponding to the unique composition $\{n+1\}$ of length 1
of $(n+1)$, we can rewrite the last expression in the following
form:
$$b_{n+1}=\sum_{c'\,\in\, {{\cal C}(n+1)}}\; \prod_{n_i\,\in\, c'}{(-a_{n_i})}\,.$$
This is Formula (\ref{enne0x}) written for the integer $(n+1)$
instead of $n$, so our proof by induction is now complete and
(\ref{enne0x}) is true for all values of $n\geq 1$.\\ Usually it
is more convenient to use partitions in place of compositions.
Using the partitions, Formula (\ref{enne0x}) becomes
\begin{eqnarray*}
 b_n&=&\sum_{p\,\in\, {{\cal P}(n)}}\,\mu(p)\; \prod_{n_i\, \in\, p}{(-a_{n_i})},\quad \forall n\geq 1\,,
\end{eqnarray*}
where the factor $\mu(p)$ is given in Equation
(\ref{numpart}).\fine\vskip 0.8cm

\section{Different approaches}

As already remarked, one may find results similar to those
expressed by Theorem \ref{uno} also in other papers. We will give
here a short historical review and with it we will stress how a
simple 
proof of Theorem \ref{uno} is not available. In fact, the result
expressed by Theorem \ref{uno} is usually seen in the context of
the theory of formal power series or in relation with the Taylor
series expansion of functions.\\

In 1825, Scherk (pages 59--66 of \cite{scherk1825}) showed how to
solve with iterated substitutions a triangular linear system like
the one corresponding to Eq. (\ref{recursion}) (that we present in
Eq. (\ref{sistscherk})), and found correctly that the $n$-th
indeterminate $b_n$ is obtained by summing up $2^{n-1}$ terms
involving products of the sequence elements $a_i$,
$\forall\,i=1,\ldots,n$. The number $2^{n-1}$ is equal to the
number of compositions of $n$, but Scherk did not
recognize the sum over the compositions of $n$ in his result.

The triangular system considered by Scherk was written in matrix
form by Brioschi in 1858 \cite{Brioschi1858}, and solved
using Cramer's rule, arriving at the
following expression for the $n$-th term of the sequence:
$$
b_n=(-1)^n\det\left(%
\begin{array}{ccccc}
  a_1 & 1 & 0 & \ldots & 0 \\
  a_2 & a_1 & 1 & \ldots & 0 \\
  a_3 & a_2 & a_1 & \ldots & 0 \\
  \vdots & \vdots & \vdots & \ddots & \vdots \\
  a_n & a_{n-1} & a_{n-2} & \ldots & a_1 \\
\end{array}%
\right)\,. $$ Brioschi then obtained the following explicit
formula:
\begin{equation}\label{andiofantea}
b_n=\sum_{q_1+2q_2+3q_3+\ldots +nq_n=n} (-1)^q\,\frac{q!}{q_1!
q_2!\cdots q_{n}!}\, {a_1^{q_1}a_2^{q_2}\cdots a_n^{q_n}}\,,
\end{equation}
where $q=q_1+q_2+\ldots +q_n$, and the sum is over all sets of
non-negative integers $\{q_1,q_2,\ldots,q_n\}$ that are solutions
of the diophantine equation
\begin{equation}\label{diofantina}
  q_1+2q_2+3q_3+\ldots +nq_n=n\,.
\end{equation}
In his work Brioschi also reported that the same formula
(\ref{andiofantea}) was obtained one year before by Fergola
\cite{Fergola1857}.

Indeed, Fergola found Eq. (\ref{andiofantea}) in 1857
\cite{Fergola1857} using an explicit formula he earlier found in
\cite{Fergola1856} for the expansion of the $n$-th derivative of
the inverse function $f^{-1}(x)$ in terms of the derivatives
of the function $f(x)$ up to the $n$-th order, where $f(x)$ is any
differentiable function. At the beginning of his paper \cite{Fergola1857}, Fergola said
that Sylvester remarked the importance of partitions in connection
with his work and the end of the same paper he stated without
proof that the number of solutions $\{q_1,q_2,\ldots,q_n\}$ of the
diophantine equation (\ref{diofantina}), with $q_i\geq 0$, is
equal to the number of partitions of $n$.

Sylvester in 1871 \cite{Sylvester1871} stated without
proof that the set of solutions of the diophantine equation
(\ref{diofantina}) are in a one to one correspondence with the set
of partitions of $n$. This fact is crucial to show that Eq.
(\ref{andiofantea}) and the result expressed by our Theorem \ref{uno} are equivalent, so
we will give a simple proof of this fact in the following
Proposition.

\begin{prop}$\!\!\!{\bf .}$\label{equivdiof}
The sets of solutions $\{q_1,q_2,\ldots,q_n\}$ of the diophantine
equation (\ref{diofantina}), with $q_i\geq 0$, are in a one to one
correspondence with the partitions $p=\{n_1,n_2,\ldots,n_k\}$ of
$n$. The correspondence is obtained by taking $q_i=m_p(n_i)$, if
$n_i\in p$, and $0$ otherwise. It then follows that
$l(p)=\sum_{i=1}^n q_i$.
\end{prop}
\noindent {\bf Proof}. Let $\{q_1,q_2,\ldots,q_n\}$ be a set of
non negative integers that are solutions of the diophantine
equation (\ref{diofantina}), then
$p=\{1q_1,2q_2,\ldots,nq_n\}\backslash\{0\}$ is a partition of
$n$. Viceversa, given a partition $p=\{n_1,n_2,\ldots\}$ of $n$,
with $n_i\geq 1$ and multiplicity $m_p(n_i)$, the set
$\{q_1,q_2,\ldots,q_n\}$, where $q_i=m_p(n_i)$, if $n_i\in p$, and
$0$ otherwise, is a set of non negative integers that are
solutions of the diophantine equation (\ref{diofantina}). By the
discussion following Eq. (\ref{numpart}) one then has that
$l(p)=\sum_{i=1}^n q_i$.\fine\\

Proposition \ref{equivdiof} can now be applied to prove the
following Proposition.

\begin{prop}$\!\!\!{\bf .}$\label{equivformula}
Eq. (\ref{andiofantea}) is equivalent to Eq. (\ref{eqt2})
.\end{prop} \noindent {\bf Proof}. Proposition
\ref{equivdiof} implies that the sum in Eq. (\ref{andiofantea}) is
over all partitions $p$ of $n$, where
$p=\{1q_1,2q_2,\ldots,nq_n\}\backslash\{0\}$. 
The factor $\frac{q!}{q_1! q_2!\cdots q_{n}!}$ in Eq.
(\ref{andiofantea}) is then equal to the factor $\mu(p)$ defined
in Eq. (\ref{numpart}), because $q=q_1+q_2+\ldots+q_n=l(p)$.
Moreover, the product $(-1)^q\,{a_1^{q_1}a_2^{q_2}\cdots
a_n^{q_n}}$ in Eq. (\ref{andiofantea}) can be rewritten as
${(-a_1)^{q_1}(-a_2)^{q_2}\cdots (-a_n)^{q_n}}$, and this is equal
to the product ${\prod_{ n_i\,\in\, p
} \left(-{a_{n_i}}\right)}$ in Eq. (\ref{eqt2}). \fine\\

Nowadays, links between the theory of formal power series and
partitions are encountered in many textbooks and papers (see, for
example, \cite{Comtet1974}). We
present, for example,  a proof of Theorem \ref{uno} in the context of the theory of formal power series.\\

\noindent {\bf Proof of Theorem \ref{uno} (formal power series
version)}. Consider the formal power series (\ref{serieformali})
satisfying $b(x)=1/a(x)$. Being $a_0=1$, one has
$a(x)=1+\sum_{i=1}^\infty a_i x^i$, so one can write:
$$b(x)=\frac{1}{a(x)}=\left(1+\sum_{i=1}^\infty a_i
x^i\right)^{-1}=1+\sum_{k=1}^\infty \left(-\sum_{i=1}^\infty a_i
x^i\right)^k
$$
where in the last member the geometric series expansion has been
used. Using the multinomial theorem
$$(x_1+x_2+\ldots)^k=\sum_{r_1+r_2+\ldots=k} \,{k \choose
{r_1,r_2,\ldots}}\; x_1^{r_1}x_2^{r_2}\cdots\,,
$$
where $r_1,r_2,\ldots$ are non-negative integers, to expand the
terms $\left(-\sum_{i=1}^\infty a_i x^i\right)^k$, and factoring
out the same powers of $x$, one may write:
$$b(x)=1+\sum_{k=1}^\infty\
\sum_{r_1+r_2+\ldots=k} \,{k \choose
{r_1,r_2,\ldots}}\; (-a_1 x)^{r_1}(-a_2 x^2)^{r_2}\cdots\,=
$$
$$=1+\sum_{k=1}^\infty\
\sum_{r_1+r_2+\ldots=k} \,x^{r_1+2r_2+\ldots}{k \choose
{r_1,r_2,\ldots}}\; (-a_1 )^{r_1}(-a_2 )^{r_2}\cdots=
$$
$$=1+\sum_{n=1}^\infty x^n \sum_{r_1+2r_2+\ldots +nr_n=n}
{r_1+r_2+\ldots +r_n \choose {r_1,r_2,\ldots,r_n}}\; (-a_1
)^{r_1}(-a_2 )^{r_2}\cdots (-a_n )^{r_n}\,.
$$
Using Proposition \ref{equivdiof} and Eq. (\ref{numpart}) one has:
$$b(x)=1+\sum_{n=1}^\infty x^n \sum_{p\,\in\, {\cal P}(n)}
\,\mu(p)\;{\prod_{  n_i\,\in\, p }(-a_{n_i})}\,,
$$
so that, for Eq. (\ref{serieformali}), one has $\forall n>0$:
$$
b_n=\sum_{p\,\in\, {\cal P}(n)} \,\mu(p)\;{\prod_{ n_i\,\in\, p
}(-a_{n_i})}\,,
$$
This expression for $b_n$ is equal to Eq. (\ref{eqt2}) in Theorem \ref{uno}.\fine\\

We listed such a different methods to obtain explicit formulas
using partitions of integers for the elements of sequences
satisfying
Condition (\ref{condition}), 
to point out that no complete elementary method, not involving
calculus or formal power series, is available in the literature.
We think moreover that our proof is simple because it only uses
some basic arithmetics and combinatorics.\\

\section{Examples}\label{examples}

We already observed that the equation
$a(x)b(x)=1$ for two formal power series is equivalent to
Condition (\ref{condition}). This applies in the particular case
of two functions $f(x)$ and $1/f(x)$ for which a Taylor series
expansion exists in $x=x_0$, and leads to following expressions for the sequence elements in terms of partitions of integers:
\begin{eqnarray}\label{entry7}
\nonumber a_n&=&\frac{1}{n!}\left.\frac{{\rm d}^n}{{\rm d}x^n}f(x)\right|_{x=x_0},\\
  b_n&=&\frac{1}{n!}\left.\frac{{\rm d}^n}{{\rm d}x^n}\frac{1}{f(x)}\right|_{x=x_0}= \sum_{p\,\in\, {\cal P}(n)} \,\mu(p)\;{\prod_{  n_i\,\in\, p}
    \left(-\frac{1}{n_i!}\left.\frac{{\rm d}^{n_i}}{{\rm d}x^{n_i}}f(x)\right|_{x=x_0}\right)}
\end{eqnarray}

Many sequences of numbers can be defined as the coefficients of the Taylor series expansions of convenient generating functions $g(x)$. In all these cases, the sequence elements can be written in terms of partitions of integers using Eq. (\ref{entry7}) for the reciprocal $f(x)=1/g(x)$ of the generating function.\\

We report in Table \ref{tab1} a few sequences of numbers
satisfying Condition (\ref{condition}). 
In the Table one should consider $n>0$, because for $n=0$ not all the listed
entries are equal to $1$, as we have supposed. The corresponding expansions in terms of partitions of integers are reported in Table \ref{tab2}. These expansions are usually found using Eq. (\ref{entry7}), but in the following we will prove some of them using a purely arithmetic method.

 {\center
\begin{table}
\begin{tabular}{|c|c|c|c|}
  \hline
 &&&\\
  Entry & Name & $a_n$ & $b_n$ \\
 &&&\\
  \hline
 &&&\\
  1& Bernoulli numbers & $\frac{1}{(n+1)!}$ & $\frac{1}{n!}\,B_n$ \\
 &&&\\
  2&Even Bernoulli numbers & $\frac{2}{(2n+2)!}$ & $-\frac{2n-1}{(2n)!}\,B_{2n}$ \\
 &&&\\
  3&Euler numbers & $\frac{1+(-1)^{n}}{2n!}$ & $\frac{1}{n!}\,E_{n}$ \\
 &&&\\
  4&Even Euler numbers & $\frac{1}{(2n)!}$ & $\frac{1}{(2n)!}\,E_{2n}$ \\
 &&&\\
  5&Fibonacci numbers & $\frac{(-1)^{n}-1}{2}$ & $F_n$ \\
 &&&\\
  6&Even Fibonacci numbers & $-n$ & $F_{2n}$ \\
  \hline
\end{tabular}\\
\caption{Some sequences satisfying Condition
(\ref{condition}).\label{tab1}}
\end{table}}
\vskip 3mm

The aim of the cited papers \cite{Fergola1857,Brioschi1858} was to
give an explicit formula for the Bernoulli and Euler numbers. For
this reason we list Entries 1--4 in Table \ref{tab1}. The proof of
Condition (\ref{condition}) for Entries 1, 3 and 4 is classical
and can easily be found in the literature (for example in
\cite{lucas1891}, starting from Eq. (3), p. 238, and Eq. (1) p.
256). We only report below a proof of Condition (\ref{condition})
for Entry 2.

Explicit formulas for the Fibonacci numbers in terms of partitions
of integers are less known. Therefore, we report below a proof of
Condition (\ref{condition}) for Entries 5 and 6.\\

\noindent {\bf Proof of Condition (\ref{condition}) for Entry 2}.
We need Condition (\ref{condition}) for Entry 1 of Table
\ref{tab1}, that is:
\begin{equation}\label{berpcond1}
\sum_{h=0}^{n}\frac{1}{(n+1-h)!}\, \frac{B_h}{h!}=0\,.
 \end{equation}
It is well-known that $B_{2n+1}=0$, $\forall n>0$. A purely
arithmetic proof of this fact is outlined in \cite{lucas1891} p.
238. It is based on a formula 
published by J. Faulhaber in 1631 and proved by B. Pascal in 1654,
as reported in \cite{beardon1996}, Theorem 3.1 (i). \footnote{A
standard proof uses the generating function for the Bernoulli numbers
$$g_B(x)=\frac{x}{e^x-1}=\sum_{h=0}^\infty B_h\, \frac{x^h}{h!}$$ and the observation that $g_B(x)+x/2$ is an
even function.}\\
Using the fact that $B_{2n+1}=0$, $\forall n>0$, it is convenient
in Eq. (\ref{berpcond1}) to separate the sums over the odd and
over the even indices. We can write Eq. (\ref{berpcond1}) for
$n=2m$ and for $n=2m+1$. When $n=2m$, Eq. (\ref{berpcond1}) gives:
$$
0=\sum_{h=0}^{2m}\frac{1}{(2m+1-h)!}\,
\frac{B_h}{h!}=
\sum_{k=0}^{m}\frac{1}{(2m+1-2k)!}\, \frac{B_{2k}}{(2k)!}\,+\,\frac{B_1}{(2m)!}\,,
$$
that is:
\begin{equation}\label{lemma1}
-\,\frac{B_1}{(2m)!}\,=\,
\sum_{k=0}^{m}\frac{2m+2-2k}{(2m+2-2k)!}\, \frac{B_{2k}}{(2k)!}
\end{equation}
When $n=2m+1$, Eq. (\ref{berpcond1}) gives:
$$
0=\sum_{h=0}^{2m+1}\frac{1}{(2m+2-h)!}\,
\frac{B_h}{h!}=
\sum_{k=0}^{m}\frac{1}{(2m+2-2k)!}\, \frac{B_{2k}}{(2k)!}\,+\,\frac{B_1}{(2m+1)!}\,,
$$
that is:
\begin{equation}\label{lemma2}
-\,\frac{B_1}{(2m)!}\,=\,\sum_{k=0}^{m}\frac{2m+1}{(2m+2-2k)!}\,
\frac{B_{2k}}{(2k)!}
\end{equation}
The difference between Eqs. (\ref{lemma2}) and (\ref{lemma1})
gives:
$$
0=\sum_{k=0}^{m}\left(\frac{2m+1}{(2m+2-2k)!}-\frac{2m+2-2k}{(2m+2-2k)!}\right)\,
\frac{B_{2k}}{(2k)!}=
$$
$$
= \sum_{k=0}^{m}\frac{2k-1}{(2m+2-2k)!}\, \frac{B_{2k}}{(2k)!}
$$
and this is equivalent to Condition (\ref{condition}) for Entry 2
of Table \ref{tab1}. (A factor 2 has been added in Entry 2 of
Table \ref{tab1} just to have $a_0=1$). \fine\vskip 0.8cm

\noindent {\bf Proof of Condition (\ref{condition}) for Entries 5 and 6}.
We first recall the standard definition of the Fibonacci numbers:
$F_1=F_2=1$, $F_n=F_{n-1}+F_{n-2}$, $\forall n>2$. We then prove
by induction that the following formulas are true $\forall n\geq
1$:
\begin{equation}\label{fibodd}
F_{2n}=\sum_{h=1}^n F_{2h-1}\,,
\end{equation}
\begin{equation}\label{fibandr}
F_{2n+1}=1+\sum_{h=1}^n F_{2h}\,,
\end{equation}
\begin{equation}\label{fibrec}
  F_{2n}=n+\sum_{h=1}^{n-1}F_{2h}(n-h)\,.
\end{equation}
Eqs. (\ref{fibodd}), (\ref{fibandr}) and (\ref{fibrec}) are true for $n=1$.
Let's suppose they are true for a certain $n\geq 2$ and let's
prove they are true also for $n$ substituted by $n+1$.
 One has:
$$F_{2(n+1)}=F_{2n+2}=F_{2n+1}+F_{2n}=F_{2n+1}+\sum_{h=1}^n F_{2h-1}=\sum_{h=1}^{n+1} F_{2h-1}\,,$$
so Eq. (\ref{fibodd}) is true for all $n$.
One has:
$$F_{2(n+1)+1}=F_{2n+3}=F_{2n+1}+F_{2n+2}=1+\sum_{h=1}^n F_{2h}+F_{2n+2}=1+\sum_{h=1}^{n+1} F_{2h}\,,$$
so Eq. (\ref{fibandr}) is true for all $n$. Using this result one
has:
$$  F_{2(n+1)}=F_{2n+2}=F_{2n+1}+F_{2n}=1+\sum_{h=1}^n
F_{2h}+n+\sum_{h=1}^{n-1}F_{2h}(n-h)=$$
$$=n+1+\sum_{h=1}^n
F_{2h}+\sum_{h=1}^{n}F_{2h}(n-h)=n+1+\sum_{h=1}^{n}F_{2h}(n+1-h),
$$
and this proves that Eq. (\ref{fibrec}) is true for all $n$.\\
Let's now form the sequence $b$, with elements $b_n=F_{2n}$, $\forall
n>0$, and $b_0=1$, and the sequence $a$, with $a_n=-n$, $\forall
n>0$, and $a_0=1$. Then Eq. (\ref{fibrec}) has the form of Eq.
(\ref{recursion}) and is clearly equivalent to Condition
(\ref{condition}) reported in Entry 6 of Table \ref{tab1}.\\
Eqs. (\ref{fibodd}) and (\ref{fibandr}) can be summarized by the following:
\begin{equation}\label{fiball}
  F_{n}=\frac{1-(-1)^n}{2}+\sum_{h=1}^{n-1}F_{h}\,\frac{1-(-1)^{n-h}}{2}\,.
\end{equation}
Let's now form the sequence $b$, with elements $b_n=F_{n}$, $\forall
n>0$, and $b_0=1$, and the sequence $a$, with $a_n=\frac{(-1)^n-1}{2}$, $\forall
n>0$, and $a_0=1$. Then Eq. (\ref{fiball}) has the form of Eq.
(\ref{recursion}) and is clearly equivalent to Condition
(\ref{condition}) reported in Entry 5 of Table \ref{tab1}.
\fine\vskip 0.8cm

 Theorem \ref{uno} can now be applied to all the
sequences listed in Table \ref{tab1} above and allows to write the explicit
formulas using
partitions reported in Table \ref{tab2}.\\

{\center
\begin{table}
\begin{tabular}{|c|c|}
  \hline
 &\\
  Entry&$b_n$ \\
 &\\
  \hline
1& $\frac{1}{n!}\,B_n=\sum_{p\,\in\, {\cal P}(n)}
\,\mu(p)\;{\prod_{n_i\,\in\, p}
\left(-\frac{1}{(n_i+1)!}\right)}$ \\
 &\\
  2&$-\frac{2n-1}{(2n)!}\,B_{2n}=\,\sum_{p\,\in\, {\cal P}(n)}
\,\mu(p)\;{\prod_{ n_i\,\in\, p  }
\left(-\frac{2}{(2n_i+2)!}\right)}$ \\
& \\
 3&$\frac{1}{n!}\,E_{n}=\sum_{p\,\in\, {\cal P}(n)}
\,\mu(p)\;{\prod_{n_i\,\in\, p}
\left(-\frac{1+(-1)^{n_i}}{2\,n_i!}\right)}$ \\
 &\\
 4&$\frac{1}{(2n)!}\,E_{2n}= \sum_{p\,\in\, {\cal P}(n)} \,\mu(p)\;{\prod_{  n_i\,\in\, p}
    \left(-\frac{1}{(2n_i)!}\right)}$ \\
 &\\
  5&$F_n= \sum_{p\,\in\, {\cal P}(n)} \,\mu(p)\;{\prod_{  n_i\,\in\, p}
    \left(-\frac{(-1)^{n_i}-1}{2}\right)}$ \\
 &\\
  6&$F_{2n}= \sum_{p\,\in\, {\cal P}(n)} \,\mu(p)\;{\prod_{  n_i\,\in\, p}
    n_i}$ \\
   \hline\end{tabular}\\
  \caption{Explicit formulas using partitions for the sequences in Table \ref{tab1}.\label{tab2}}
  \end{table}}
  \vskip 3mm

Finally, we present two explicit calculations using ${\cal P}(4)$. One
has:
$${\cal P}(4)=\{\{4\},\,\{3,1\},\,\{2,2\},\,\{2,1,1\},\,\{1,1,1,1\} \}\,,$$
so from the first and sixth rows of Table \ref{tab2} one finds:
$$B_4=4!\,\left[ \frac{1!}{1!}\cdot \frac{-1}{5!}+\frac{2!}{1!\, 1!}\cdot \frac{(-1)^2}{4!\, 2!}+ \frac{2!}{2!}\cdot \frac{(-1)^2}{3!\, 3!}+ \frac{3!}{1!\, 2!}\cdot \frac{(-1)^3}{3!\,2!\,2!}+\frac{4!}{4!}\cdot \frac{(-1)^4}{2!\,2!\,2!\,2!}\right]=-\frac{1}{30}\,,$$
$$F_8=\,\left[ \frac{1!}{1!}\cdot 4+\frac{2!}{1!\, 1!}\cdot 3\cdot 1+ \frac{2!}{2!}\cdot 2\cdot 2+ \frac{3!}{1!\, 2!}\cdot 2\cdot 1\cdot 1+\frac{4!}{4!}\cdot 1\cdot 1\cdot 1\cdot 1\right]=21\,.$$

\end{document}